\documentclass[wcp]{jmlr}
\usepackage{amsmath,amssymb,bbm,algorithm,algorithmic,graphicx}

\usepackage{booktabs}

\jmlrvolume{vol}
\jmlryear{2011}
\jmlrworkshop{25th Annual Conference on Learning Theory}

\title[KL-UCB: bounded bandits, and beyond]{The KL-UCB Algorithm for Bounded Stochastic Bandits and Beyond}

  \author{\Name{Aur\'elien Garivier} \and
   \Name{Olivier Capp\'e} \Email{garivier,cappe@telecom-paristech.fr}\\
   \addr LTCI, CNRS \& Telecom ParisTech, Paris, France}

\editor{Sham Kakade, Ulrike von Luxburg}

\newcommand{\E}{\mathbb{E}} 
\newcommand{\Var}{\mathbb{V}ar} 
\newcommand{\N}{\mathbb{N}}
\renewcommand{\P}{\mathbb{P}} 
\newcommand{\R}{\mathbb{R}}

\newcommand{\1}{\mathbbm{1}}

\def\gets{\leftarrow}
\def\arm{\hbox{arm}}
\def\reward{\hbox{REWARD}}
\def\argmax{\operatorname*{arg\,max}}

\def\F{\mathcal{F}}

\begin{document}

\maketitle

\begin{abstract}
  This paper presents a finite-time analysis of the KL-UCB algorithm, an online, horizon-free
  index policy for stochastic bandit problems.
  We prove two distinct results: first, for arbitrary bounded rewards, the KL-UCB algorithm
  satisfies a uniformly better regret bound than UCB and its variants; second, in the special case of
  Bernoulli rewards, it reaches the lower bound of Lai and Robbins.
  Furthermore, we show that simple adaptations of the KL-UCB algorithm are also optimal for
  specific classes of (possibly unbounded) rewards, including those generated from exponential
  families of distributions.
  A large-scale numerical study comparing KL-UCB with its main competitors (UCB, MOSS,
  UCB-Tuned, UCB-V, DMED) shows that KL-UCB is remarkably efficient and stable, including for short time horizons. KL-UCB is also the only method that always performs better
  than the basic UCB policy.
  Our regret bounds rely on deviations results of independent interest which are stated and proved
  in the Appendix. As a by-product, we also obtain an improved regret bound for the standard UCB
  algorithm.
\end{abstract}

\section{Introduction}\label{sec:intro}

The multi-armed bandit problem is a simple, archetypal setting of reinforcement
learning, where an agent facing a slot machine with several arms tries to
maximize her profit by a proper choice of arm draws.  In the
stochastic\footnote{Another interesting setting is the \emph{adversarial}
  bandit problem, where the rewards are not stochastic but chosen by an
  opponent - this setting is not the subject of this paper.} bandit problem,
the agent sequentially chooses, for $t=1,2,\dots, n$, an arm
$A_t\in\{1,\dots,K\}$, and receives a reward $X_t$ such that, conditionally on
the arm choices $A_1,A_2,\dots$, the rewards are independent and identically
distributed, with expectation $\mu_{A_1},\mu_{A_2},\dots$. Her \emph{policy} is
the (possibly randomized) decision rule that, to every past observations
$(A_{1}, X_1,\dots,A_{t-1}, X_{t-1})$, associates her next choice $A_t$. The best
choice is any arm $a^*$ with maximal expected reward $\mu_{a^*}$. The
performance of her policy can be measured by the \emph{regret} $R_n$, defined as the
difference between the rewards she accumulates up to the horizon $t=n$, and the rewards that she would have accumulated during the same
period, had she known from the beginning which arm had the highest expected reward. 

The agent typically faces an ``exploration versus exploitation dilemma''~: at
time $t$, she can take advantage of the information she has gathered, by
choosing the so-far best performing arm, but she has to consider the
possibility that the other arms are actually under-rated and she must play
sufficiently often all of them.  Since the works of~\cite{Gittins79bandits} in
the 1970s, this problem raised much interest and several variants, solutions
and extensions have been proposed, see \cite{EvendarMannorMansour02PACbandits} and references therein.

Two families of bandit settings can be distinguished: in the first family, the
distribution of $X_t$ given $A_t=a$ is assumed to belong to a family
$\{p_\theta, \theta \in \Theta_a\}$ of probability distributions.  In a
particular parametric framework, \cite{LaiRobbins85bandits} proved a
lower-bound on the performance of any policy, and determined optimal policies.
This framework was extended to multi-parameter models
by \cite{BurnetasKatehakis97bandits}  who showed that the number
of draws up to time $n$, $N_a(n)$, of any sub-optimal arm $a$ is lower-bounded by
\begin{equation}\label{eq:binfNb}
  N_a(n) \geq \left( \frac{1}{\inf_{\theta \in \Theta_a : E[p_\theta]>\mu_{a^*} } KL\left( p_{\theta_a}, p_{\theta}\right)} +o(1)\right)\;\log(n),
\end{equation}
where $KL$ denotes the Kullback-Leibler divergence and $E(p_\theta)$ is
the expectation under $p_\theta$; hence, the regret is lower-bounded as
follows:
\begin{equation}\label{eq:binfRegret}\liminf_{n\to\infty}
  \frac{\E[R_n]}{\log(n)} \geq \sum_{a:\mu_a<\mu_{a^*}} 
  \frac{\mu_{a^*}-\mu_a}{\inf_{\theta \in \Theta_a : E[p_\theta]>\mu_{a^*} }KL\left( p_{\theta_a}, p_{\theta}\right)}\;.
\end{equation}
Recently, \cite{HondaTakemura10DMED} proposed an algorithm called
\emph{Deterministic Minimum Empirical Divergence (DMED)} which they proved to be first order
optimal. This algorithm, which maintains a list of arms that are close enough
to the best one (and which thus must be played), is inspired by large
deviations ideas and relies on the availability of the rate function associated to the reward distribution.

In the second family of bandit problems, the rewards are only
assumed to be bounded (say, between $0$ and $1$), and policies rely directly on
the estimates of the expected rewards for each arm. The KL-UCB algorithm considered in this paper is primarily meant to address this second, non-parametric, setting.
We will nonetheless show that KL-UCB also matches the lower bound of \cite{BurnetasKatehakis97bandits} in the binary case and that it can be extended to a larger class of parametric bandit problems.

Among all the bandit policies that have been proposed, a particular family has
raised a strong interest, after \cite{Gittins79bandits} proved that (in
the Bayesian setting he considered) optimal policies could be chosen in the
following very special form: compute for each arm a \emph{dynamic allocation
  index} (which only depends on the draws of this arm), and simply choose an
arm with maximal index. These policies not only compute an estimate of the
expected rewards, but rather an \emph{upper-confidence bound} (UCB), and the
agent's choice is an arm with highest UCB. This approach is sometimes called
``optimistic'', as the agent acts as if, at each instant, the expected rewards 
were equal to the highest possible values that are compatible with her past observations.  \cite{AuerEtAl02FiniteTime}, following \cite{Agrawal:95},
proposed and analyzed two variants, UCB1 (usually called simply UCB in latter works) and UCB2, for which they provided regret
bounds. UCB is an online, horizon-free procedure for which \citep{AuerEtAl02FiniteTime} proves that there exists a constant $C$ such that
\begin{equation}\label{eq:UCBAuer}
\E[R_n] \leq \sum_{a:\mu_a<\mu_{a^*}} \frac{8\log(n)}{\left( \mu_{a^*}-\mu_{a} \right)} + C \; .
\end{equation}
The UCB2 variant relies on a parameter $\alpha$ that needs to
be tuned, depending in particular on the horizon, and satisfies the tighter regret bound
\begin{equation*}
\E[R_n] \leq \sum_{a:\mu_a<\mu_{a^*}} \frac{(1+\epsilon(\alpha))\log(n)}{2\left(\mu_{a^*}-\mu_{a} \right)} + C(\alpha) \; ,
\end{equation*}
where $\epsilon(\alpha) > 0$ is a constant that can get arbitrary small when
$\alpha$ is small, at the expense of an increased value of the constant $C(\alpha)$.
The constant $1/2$ in front of the factor $\log(n)/\left(\mu_{a^*}-\mu_{a}
\right)$ cannot be improved. We show in Proposition~\ref{prop:betterUCB}, as a by-product of our analysis, that a correctly tuned UCB algorithm satisfies a similar bound.
However, \cite{AuerEtAl02FiniteTime} found in numerical experiments that UCB and UCB2 were outperformed by a third heuristic variant called UCB-Tuned, which includes estimates of the variance, but no theoretical guarantee was given. In a latter work, \cite{AudibertEtAlUCBV} proposed a related policy, called UCB-V, which uses an empirical version of the Bernstein bound to obtain refined upper confidence bounds. Recently, \cite{audibert:bubeck:2010:minimax} introduced an improved UCB algorithm, termed MOSS, which achieves the distribution-free optimal rate.

In this contribution, we first consider the stochastic, non-parametric, bounded bandit problem. We consider an online index policy, called KL-UCB (for Kullback-Leibler UCB), that requires no problem- or horizon-dependent tuning. This algorithm was recently advocated by \cite{filippi10}, together with a similar procedure for Markov Decision Processes~\citep{allerton10}, and we learnt since our initial submission that an analysis of the Bernoulli case can also be found in~\cite{MaillardMunosStoltz11klucb}, together with other results. We prove in Theorem~\ref{th:borneRegret} below
that the regret of KL-UCB satisfies
\[
 \limsup_{n\to\infty} \frac{\E[R_n]}{\log(n)} \leq \sum_{a:\mu_a<\mu_{a^*}}
\frac{\mu_{a^*}-\mu_a}{d(\mu_a, \mu_{a^*})}\;,
\]
where $d(p,q) = p\log(p/q)+(1-p)\log((1-p/(1-q))$ denotes the Kullback-Leibler
divergence between Bernoulli distributions of parameters $p$ and $q$, respectively. This comes as a consequence of Theorem~\ref{th:borneNbplayed}, a
non-asymptotic upper-bound on the number of draws of a sub-optimal arm $a$: for
all $\epsilon>0$ there exist $C_1, C_2(\epsilon)$ and $\beta(\epsilon)$ such
that
\[
 \E[N_n(a)] \leq \frac{\log(n)}{d(\mu_a, \mu_{a^*})}(1+\epsilon) + C_1
\log(\log(n)) + \frac{C_2(\epsilon)}{n^{\beta(\epsilon)}}\;.
\]
We insist on the
fact that, despite the presence of divergence $d$, this bound is not specific to
the Bernoulli case and applies to all reward distributions bounded in $[0,1]$
(and thus, by rescaling, to all bounded reward distributions). By Pinsker's
inequality, $d(\mu_a,\mu_{a^*}) > 2 (\mu_a-\mu_{a^*})^2$, and thus KL-UCB has strictly better theoretical guarantees than UCB, while it has the same range of application. The improvement appears to be significant in simulations. Moreover, KL-UCB is the first index policy that reaches the lower-bound of
\cite{LaiRobbins85bandits} for binary rewards;  it does also achieve lower 
regret than UCB-V in that case. Hence, KL-UCB is both a general-purpose
procedure for bounded bandits, and an optimal solution for the binary case.

Furthermore, it is easy to adapt KL-UCB to particular (possibly non-bounded) bandit settings, when the distribution of reward is known to belong to some family of probability laws. By simply changing the definition of the divergence $d$, optimal algorithms may be built in a great variety of situations.

The proofs we give for these results are short and elementary. They rely on
deviation results for bounded variables that are of independent interest~:
Lemma~\ref{lem:bounded2bernoulli} shows that Bernoulli variable are, in a
sense, the ``least concentrated'' bounded variables with a given expectation (as is
well-known for variance), and Theorem~\ref{th:borneDev} shows an efficient way
to build confidence intervals for sums of bounded variables with an unknown
number of summands. 

In practice, numerical experiments confirm the significant advantage of KL-UCB
over existing procedures; not only does this method outperform UCB, MOSS,
UCB-V and even UCB-tuned in various scenarios, but it also compares well to DMED in the Bernoulli case, especially for small or moderate horizons.

The paper is organized as follows: in Section~\ref{sec:algo}, we introduce some
notation and present the KL-UCB algorithm.  Section~\ref{sec:bound} contains the main results of the paper, namely the regret bound for KL-UCB and
the optimality in the Bernoulli case. In Section~\ref{sec:parametricKLUCB}, we show how to adapt the KL-UCB algorithm to address general families of reward distributions, and we provide finite-time regret bounds showing asymptotic optimality. Section~\ref{sec:simu} reports the
results of extensive numerical experiments, showing the practical superiority of
KL-UCB. Section~\ref{sec:proof} is devoted to an elementary proof of the main
theorem. 
Finally, the Appendix gathers some deviation results that are useful in the proofs of our regret bounds, but which are also of independent interest.

\section{The KL-UCB Algorithm}\label{sec:algo}
We consider the following bandit problem: the set of actions is
$\{1,\dots,K\}$, where $K$ denotes a finite integer.  For each
$a\in\{1,\dots,K\}$, the rewards $(X_{a,t})_{t\geq 1}$ are independent and
bounded\footnote{if the rewards are bounded in another interval $[a,b]$, they should first be rescaled to $[0,1]$.} in $\Theta = [0,1]$ with common expectation $\mu_a$. The sequences
$(X_{a,\cdot})_a$ are independent.  At each time step $t=1,2,\dots$, the agent
chooses an action $A_t$ according to his past observations (possibly using some
independent randomization) and gets the reward $X_t = X_{A_t, N_{A_t}(t)}$,
where $N_a(t) = \sum_{s=1}^t \1\{A_s = a\}$ denotes the number of times action $a$ was chosen up to time $t$. The sum of rewards she has obtained when choosing action $a$ is denoted by $S_a(t) = \sum_{s\leq t} \1\{A_s = a\} X_s =
\sum_{s=1}^{N_a(t)} X_{a,s}$. 
For $(p,q)\in \Theta^2$ denote the Bernoulli
Kullback-Leibler divergence by
\[d(p,q) = p\log\frac{p}{q} + (1-p)\log\frac{1-p}{1-q}\;,\] with, by convention, $0\log0 = 0\log0/0 =0$ and
$x\log x/0=+\infty$ for $x>0$.

\begin{algorithm}
  \caption{KL-UCB}
  \label{alg:klUCB}
  \begin{algorithmic}[1]
    \REQUIRE $n$ (horizon), $K$ (number of arms), $\reward$ (reward function,
    bounded in $[0,1]$) \FOR {$t = 1$ \textbf{to} $K$} \STATE $N[t] \gets 1$
    \STATE $S[t] \gets \reward(\arm=t)$
    \ENDFOR
    \FOR {$t = K+1$ \textbf{to} $n$} \STATE $a \gets \argmax_{1\leq a\leq K}
    \max\left\{q\in\Theta : N[a]\,d\left( \frac{S[a]}{N[a]}, q \right) \leq
      \log(t) + c\log(\log(t))\right\}$ \STATE $r \gets \reward(\arm=a)$ \STATE
    $N[a] \gets N[a]+1$ \STATE $S[a] \gets S[a]+r$
    \ENDFOR
  \end{algorithmic}
\end{algorithm}

Algorithm~\ref{alg:klUCB} provides the pseudo-code for KL-UCB. On line 6, $c$ is a parameter that, in the regret
bound stated below in Theorem~\ref{th:borneRegret} is chosen equal to $3$; in
practice, however, we recommend to take $c=0$ for optimal performance.
For each arm $a$ the upper-confidence bound
\[\max\left\{q\in\Theta : N[a]\,d\left( \frac{S[a]}{N[a]}, q \right) \leq
  \log(t) + c\log(\log(t))\right\}\] 
can be efficiently computed using Newton iterations, as for any $p\in[0,1]$ the
function $q\mapsto d(p,q)$ is strictly convex and increasing on the
interval $[p, 1]$. In case of ties between severals arms, any maximizer can be
chosen (for instance, at random).
The KL-UCB elaborates on ideas suggested in Sections~3 and~4 of \cite{LaiRobbins85bandits}.

\section{Regret bounds and optimality} \label{sec:bound} We first state the
main result of this paper.  It is a direct consequence of the non-asymptotic
bound in Theorem~\ref{th:borneNbplayed} stated below.
\begin{theorem}\label{th:borneRegret}
  Consider a bandit problem with $K$ arms and independent rewards bounded in
  $[0,1]$, and denote by $a^*$ an optimal arm.  Choosing $c=3$, the regret of
  the KL-UCB algorithm satisfies:
  \[\limsup_{n\to\infty} \frac{\E[R_n]}{\log(n)} \leq \sum_{a:\mu_a<\mu_{a^*}}
  \frac{\mu_{a^*}-\mu_a}{d(\mu_a, \mu_{a^*})}\;.\]
\end{theorem}

\begin{theorem}\label{th:borneNbplayed}
  Consider a bandit problem with $K$ arms and independent rewards  bounded in
  $[0,1]$. Let $\epsilon>0$, and take $c=3$ in Algorithm~\ref{alg:klUCB}.  Let $a^*$
  denote an arm with maximal expected reward $\mu_{a^*}$, and let $a$ be an arm
  such that $\mu_a < \mu_{a^*}$.  For any positive integer $n$, the number of
  times algorithm KL-UCB chooses arm $a$ is upper-bounded by
  \[\E[N_n(a)] \leq \frac{\log(n)}{d(\mu_a, \mu_{a^*})}(1+\epsilon) + C_1
  \log(\log(n)) + \frac{C_2(\epsilon)}{n^{\beta(\epsilon)}}\;,\] where $C_1$
  denotes a positive constant and where $C_2(\epsilon)$ and $\beta(\epsilon)$
  denote positive functions of $\epsilon$. Hence,
  \[\limsup_{n\to\infty} \frac{\E[N_n(a)]}{\log(n)} \leq \frac{1}{d(\mu_a,
    \mu_{a^*})}\;.\]
\end{theorem}

\begin{corollary}
  If the reward distributions are Bernoulli, the KL-UCB algorithm is
  asymptotically optimal, in the sense that the regret of KL-UCB matches the
  lower-bound proved by~\cite{LaiRobbins85bandits} and generalized
  by~\cite{BurnetasKatehakis97bandits}:
  \[N_n(a) \geq \left( \frac{1}{d(\mu_a, \mu_{a^*})} + o(1) \right)\log(n)\] with a
  probability tending to $1$.
\end{corollary}

The KL-UCB algorithm thus appears to be (asymptotically) optimal for Bernoulli rewards. However, Lemma~\ref{lem:bounded2bernoulli} shows that the Chernoff bounds obtained for Bernoulli variables actually apply to any variable with range $[0,1]$. This is why 
KL-UCB is not only efficient in the binary case, but also for general bounded rewards.

As a by-product, we obtain an improved upper-bound for the regret of the  UCB algorithm:
\begin{proposition}\label{prop:betterUCB}
Consider the UCB algorithm tuned as follows: at step $t>K$, the arm that maximizes the upper-bound $S[a]/N[a] + \sqrt{(\log(t)+c\log\log(t))/(2 N[a])}$ is chosen. 
Then, for the choice $c=3$, the number of draws of a sub-optimal arm $a$ is upper-bounded as:
  \begin{equation}\label{eq:UCBmeilleureborne}
   \E[N_n(a)] \leq \frac{\log(n)}{2(\mu_a-\mu_{a^*})^2}(1+\epsilon) + C_1
  \log(\log(n)) + \frac{C_2(\epsilon)}{n^{\beta(\epsilon)}}\;.
\end{equation}
\end{proposition}
This bound is ``optimal'', in the sense that the constant $1/2$ in the logarithmic term cannot be improved. 
The proof of this proposition just mimics that of Section~\ref{sec:proof} (which concerns KL-UCB), using the quadratic divergence  $d(p,q):= 2(p-q)^2$ instead of the Kullback-Leibler divergence; it is thus omitted.
In contrast, Pinsker's inequality $d(\mu_a, \mu_{a^*})\geq 2(\mu_a-\mu_{a^*})^2$ shows that KL-UCB dominates UCB, and we will see in the simulation study that the difference is significant, even for smaller values of the horizon.

\begin{remark}
  At line $6$, Algorithm~\ref{alg:klUCB} computes for each arm
  $a\in\{1,\dots,K\}$ the upper-confidence bound
  \[\max\left\{q\in\Theta : N[a] \, d\left( \frac{S[a]}{N[a]}, q \right)
    \leq \log(t) + c\log(\log(t))\right\}\;.\]
  The level of this confidence bound is parameterized by the exploration function $\log(t)+c\log(\log(t))$, and the
  results of Theorems~\ref{th:borneRegret} and \ref{th:borneNbplayed} are true
  as soon as $c\geq3$.  However, similar results can be proved with an exploration function equal to  $(1+\epsilon)\log(t)$ (instead of
  $\log(t)+c\log(\log(t))$) for every $\epsilon>0$; this is no surprise, as
  $(1+\epsilon)\log(t) \geq \log(t)+c\log(\log(t))$ when $t$ is large
  enough. But ``large enough'', in that case, can be quite large~: for
  $\epsilon=0.1$, this holds true only for $t > 2.10^{51}$. This is why, in
  practice (and for the simulations presented in
  Section~\ref{sec:simu}), we rather suggest to choose $c=0$. 
\end{remark}

\section{KL-UCB for parametric families of reward distributions}\label{sec:parametricKLUCB}
The KL-UCB algorithm makes no assumption on the distribution of the rewards,
except that they are bounded. Actually, the definition of the divergence
function $d$ in KL-UCB is dictated by the rate function of the Large Deviations
Principle satisfied by Bernoulli random variables: the proof of
Theorem~\ref{th:borneDev} relies on the control of the Fenchel-Legendre
transform of the Bernoulli distribution. Thanks to
Lemma~\ref{lem:bounded2bernoulli}, this choice also makes sense for all
bounded variables.

But the method presented here is not limited to the Bernoulli case:
 KL-UCB can very easily be adapted to other settings by
choosing an appropriate divergence function $d$.
As an illustration, we will assume in this section that, for each arm $a$, the distribution of rewards belongs to
a \emph{canonical exponential family}, i.e., that its density with respect to some
reference measure can be written as $p_{\theta_a}(x)$ for some real parameter
$\theta_a$, with
\begin{equation}\label{eq:densiteModeleExpo}p_{\theta}(x) = \exp\left( x\theta -
    b(\theta) + c(x) \right)\;,\end{equation}
where $\theta$ is a real parameter, $c$ is a real function and the $\log$-partition function $b(\cdot)$ is assumed to be twice differentiable. This family contains for instance the Exponential, Poisson, Gaussian (with fixed variance), Gamma (with fixed shape parameter) distributions (as well as, of course, the Bernoulli distribution).
For a random variable $X$ with density defined in $\eqref{eq:densiteModeleExpo}$, it is easily checked that $\mu(\theta)\triangleq \E_\theta[X] =\dot{b}(\theta)$; moreover, as $\ddot{b}(\theta)=\Var(X)>0$, the function $\theta\mapsto \mu(\theta)$ is one-to-one.
Theorem~{\ref{th:borneDevGen} (in the Appendix) states that the probability of under-estimating the performance of the best arm can be upper-bounded just as in the Bernoulli case by replacing the divergence $d(\cdot,\cdot)$ in line 6 of the KL-UCB algorithm by
\[ d(x, \mu(\theta)) = \sup_{\lambda} \left\{ \lambda x - \log \E_\theta\left[\exp(\lambda X)\right] \right\} \;.\]
For example, in the case of exponential rewards, one should choose $d(x,y) =x/y - 1 - \log (x/y)$. Or, for Poisson rewards, the right choice is $d(x,y)=y-x+x\log(x/y)$.
Then, all the results stated above apply (as the proofs do not involve the particular form of the function $d$), and in particular~:
  \[\limsup_{n\to\infty} \frac{\E[R_n]}{\log(n)} \leq \sum_{a:\mu_a<\mu_{a^*}}
  \frac{\mu_{a^*}-\mu_a}{d(\mu_a, \mu_{a^*})}\;.\]
In order to prove that, for those families of rewards, this version of the KL-UCB algorithm matches the bound of \cite{LaiRobbins85bandits} , it remains only to show that $d(x, y) = KL(p_{\mu^{-1}(x)}, p_{\mu^{-1}(y)})$. This is the object of
Lemma~\ref{lem:dIsKl}. 
Generalizations to other families of reward distributions (possibly different from arm to arm) are conceivable, but require more technical, topological discussions, as in \cite{BurnetasKatehakis97bandits} and \cite{HondaTakemura10DMED}.

To conclude, observe that it is not required to work with the divergence function $d$ corresponding exactly to the family of reward distributions: using an upper-bound instead often leads to more simple and versatile policies at the price of only a slight loss in performance. This is illustrated in the third scenario of the simulation study presented in Section~\ref{sec:simu}, but also by Theorems~\ref{th:borneRegret} and \ref{th:borneNbplayed} for bounded variables.

\begin{lemma}\label{lem:dIsKl}
  Let $(\beta,\theta)$ be two real numbers, let $p_\beta$ and $p_\theta$ be two
  probability densities of the canonical exponential family defined
  in~\eqref{eq:densiteModeleExpo}, and let $X$ have density $p_\theta$.  Then
  Kullback-Leibler divergence $KL(p_\beta, p_\theta)$ is equal to  
  $d(\mu(\beta), \mu(\theta))$. 
  More precisely, 
  \[KL(p_\beta, p_\theta) = d(\mu(\beta), \mu(\theta))  = \mu(\beta)\left( \beta-\theta \right) -b(\beta)+b(\theta)\;.\]
\end{lemma}
\begin{proof}
  First, it holds that
  \begin{align*}
   KL(p_\beta, p_\theta)& = \int \exp\left( x\beta - b(\beta) + c(x) \right)
  \left\{ x\left( \beta-\theta \right) -b(\beta)+b(\theta) \right\} dx \\ 
  &=
  \mu(\beta)\left( \beta-\theta \right) -b(\beta)+b(\theta)\;.   
  \end{align*}
   Then, observe
  that
 $\E\left[ \exp(\lambda X) \right] = \int \exp\left( x(\beta+\lambda) -
    b(\beta) + c(x) \right) dx = \exp( b(\beta+\lambda) - b(\beta))\;.$ 
  Thus,
  for every $x$ the maximum of the (smooth, concave) function
  \[\lambda\mapsto\lambda x - \log\E\left[ \exp(\lambda X) \right] = \lambda x -
  b(\theta+\lambda) + b(\theta)\] is reached for $\lambda=\lambda^*$ such that  $x = \dot{b}(\theta +
  \lambda^*) = \mu(\theta+\lambda^*)$.  Thus, if $x=\mu(\beta)$, the fact that
  $\mu$ is one-to-one implies that $\theta+\lambda^*=\beta$ and thus that:
  \[d(\mu(\beta), \mu(\theta)) = \left( \beta-\theta \right) \mu(\beta) -
  b(\beta) + b(\theta) \;.\]
\end{proof}

\section{Numerical experiments and comparisons of the policies}\label{sec:simu}

Simulations studies 
require particular attention in the case of bandit algorithms.  As pointed out
by \cite{AudibertEtAlUCBV}, for a fixed horizon $n$ the distribution of the
regret is very poorly concentrated around its expectation. This can be
explained as follows: most of the time, the estimates of all arms remain
correctly ordered for almost all instants $t=1,\dots,n$ and the regret is of
order $\log(n)$. But sometimes, at the beginning, the best arm is
under-estimated while one of the sub-optimal arms is over-estimated, so that the
agent keeps playing the latter; and as she neglects the best arm, she has
hardly an occasion to realize her mistake, and the error perpetuates for a
very long time. This happens with a small, but not negligible probability, because
the regret is very important (of order $n$) on these occasions. Bandit
algorithms are typically designed to control the probability of such adverse
events but usually at a rate which only decays slightly faster than $1/n$, which results in very skewed regret distributions with slowly vanishing upper tails.

\begin{figure}[hbt]
  \centering
  \begin{minipage}[t]{0.44\textwidth}
  \vspace*{0pt}
  \includegraphics[height=0.8\linewidth]{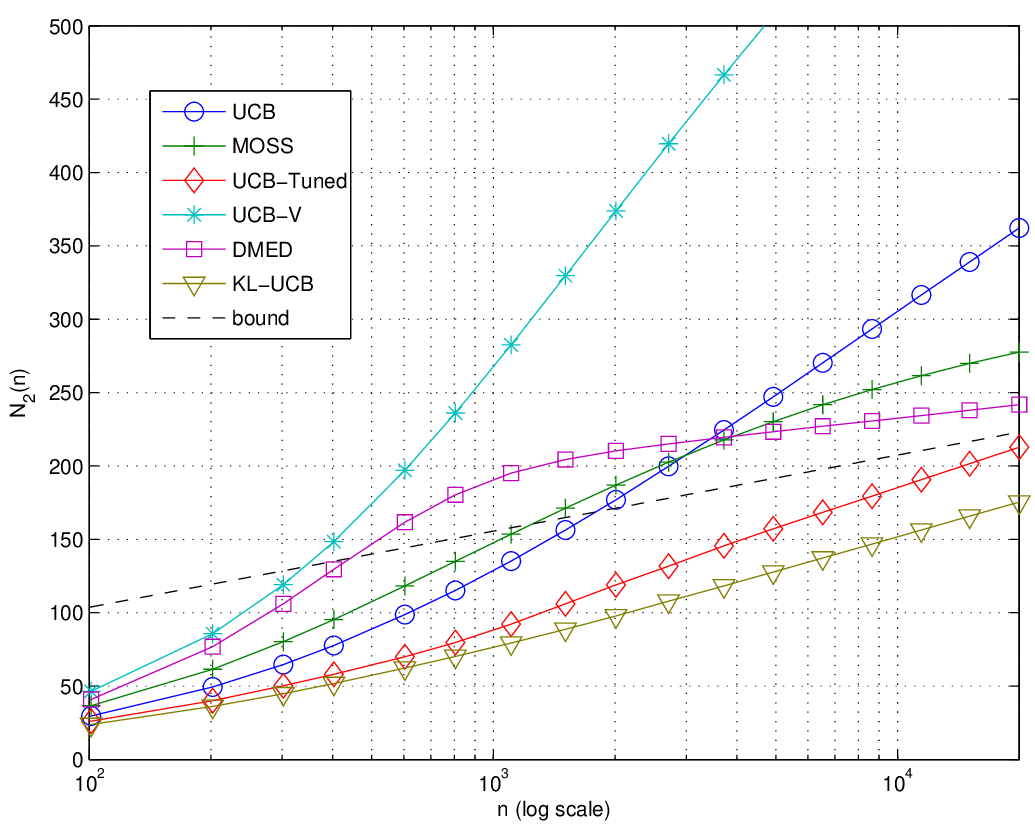}
  \end{minipage} \qquad
  \begin{minipage}[t]{0.44\textwidth}
    \vspace*{0pt}
   \includegraphics[height=0.75\linewidth]{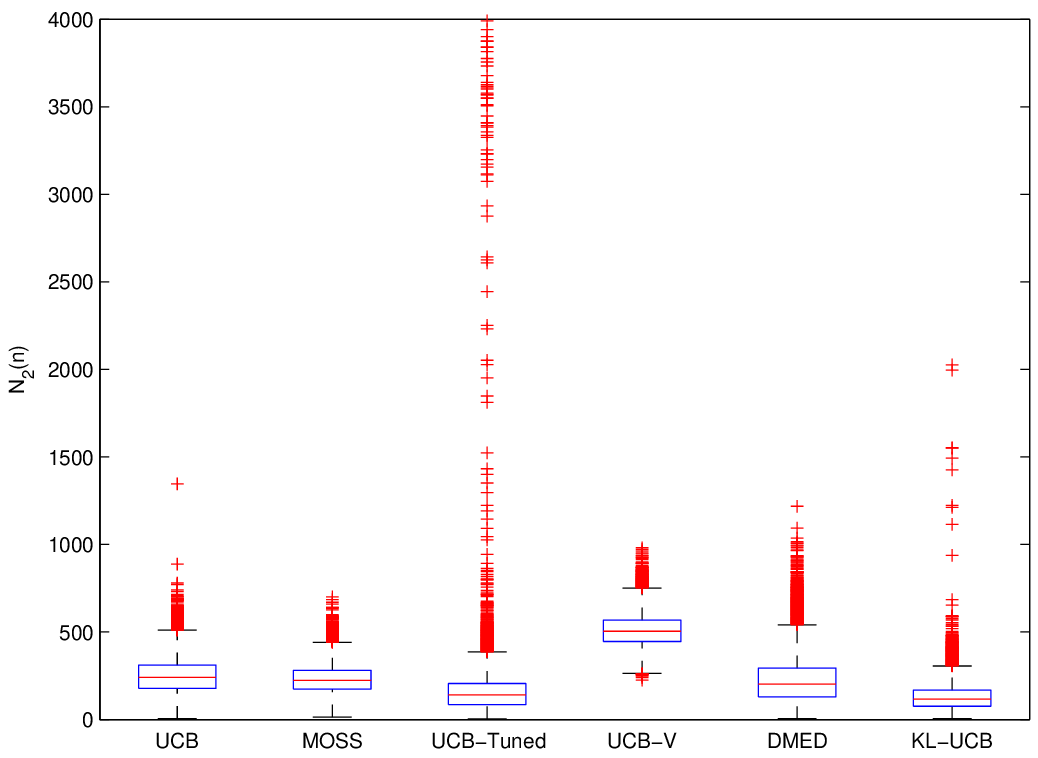}
  \end{minipage}
  \caption{Performance of the various algorithms in the simple two arms, scenario. Left, mean number of draws of the suboptimal arm as a function of time; right, box-plots showing the distributions of the number of draws of the suboptimal arm at time $n=5,000$. Results based on $50,000$ independent runs.}
  \label{fig:2a}
\end{figure}

\subsection{Scenario 1: two arms}
We first consider the basic two arm scenario with Bernoulli rewards of
expectations $\mu_1=0.9$ and $\mu_2=0.8$, respectively. The left panel of
Figure~\ref{fig:2a} shows the average number of draws of the suboptimal arm as
a function of time (on a logarithmic scale) for KL-UCB compared to five
benchmark algorithms (UCB, MOSS, UCB-Tuned, UCB-V and DMED). The right panel of
Figure~\ref{fig:2a} shows the empirical distributions of suboptimal draws,
represented as box-and-whiskers plots, at a particular time ($t=5,000$) for
all six algorithms. These plot are obtained from $N=50,000$ independent runs of
the algorithms and the right panel of Figure~\ref{fig:2a} clearly highlight
the tail effect mentioned above. On this very simple example, we observed that
results obtained from $N=1,000$ or less simulations were not reliable,
typically resulting in a significant over-estimation\footnote{Incidentally, Theorem~\ref{th:borneDev} could be used to construct sharp confidence bounds for the regret.} of the performance of
``risky'' algorithms, in particular of UCB-Tuned. More generally, results
obtained in configurations where $N$ is much smaller than $n$ are likely to be
unreliable. For this reason, we limit our investigations to a final instant of
$n=20,000$. Note however that the average number of suboptimal draws of most
algorithms at $n=20,000$ is only of the order of $300$, showing that there is no
point in considering larger horizons for such a simple problem.

MOSS, UCB-Tuned and UCB-V are run exactly as described by \cite{audibert:bubeck:2010:minimax}, \cite{AuerEtAl02FiniteTime}
and \cite{AudibertEtAlUCBV}, respectively. For UCB, we use an upper confidence bound 
$S[a]/N[a] + \sqrt{\log(t)/(2 N[a])}$, as in Proposition~\ref{prop:betterUCB}, again with $c=0$.
Note that in our two arm scenario, $\{2(\mu_1-\mu_2)^2\}^{-1} = 50$ while
$d^{-1}(\mu_2,\mu_1) = 22.5$. Hence, the performance of DMED and KL-UCB should be
about two times better than that of UCB. The left panel of
Figure~\ref{fig:2a} does show the expected behavior but with a difference
of lesser magnitude. Indeed, one can observe that the bound
$d^{-1}(\mu_2,\mu_1)\log(n)$ (shown in dashed line) is quite pessimistic for
the values of the horizon $n$ considered here as the actual performance of
KL-UCB is significantly below the bound. 
For DMED, we follow \cite{HondaTakemura10DMED} but using
\begin{equation}
  \label{eq:dmed}
  N[a] \, d\left(\frac{S[a]}{N[a]}, \max_b\frac{S[b]}{N[b]}\right) < \log t   
\end{equation}
as the criterion to decide whether arm $a$ should be included in the list of
arms to be played. This criterion is clearly related to the decision rule used
by KL-UCB when $c=0$ (see line 6 of Algorithm~\ref{alg:klUCB}) except for the
fact that in KL-UCB the estimate $S[a]/N[a]$ is not compared to that of the
current best arm $\max_b S[b]/N[b]$ but to the corresponding upper confidence
bound. As a consequence, any arm that is not included in the list of arms to be
played by DMED would not be played by KL-UCB either (assuming that both
algorithms share a common history). As one can observe on the left panel of
Figure~\ref{fig:2a}, this results in a degraded performance for DMED. We also
observed this effect on UCB, for instance, and it seems that index algorithms
are generally preferable to their ``arm elimination'' variant.

The original proposal of \cite{HondaTakemura10DMED} consists in using in the
exploration function a factor $\log(t/N[a])$ instead of $\log(t)$, as in the MOSS algorithm.
As will be seen below on Figure~\ref{fig:10}, this variant (which we refer to as
DMED+) indeed outperforms DMED. But our previous conjecture appears to hold also
in this case as the heuristic variant of KL-UCB (termed KL-UCB+) in which
$\log(t)$ in line 6 of Algorithm~\ref{alg:klUCB} is replaced by $\log(t/N[a])$
remains preferable to DMED+. 

As final comments on Figure~\ref{fig:2a}, first
note that UCB-Tuned performs as expected ---though slightly worse than KL-UCB---
but is a very risky algorithm: the right panel of Figure~\ref{fig:2a} 
casts some doubts on the fact that the tails of $N_a(n)$ are indeed controlled
uniformly in $n$ for UCB-Tuned. 
Second, the performance of UCB-V is somewhat disappointing. Indeed, the upper-confidence bounds of UCB-V differ from those of UCB-Tuned simply by the non-asymptotic correction term $3\log(t)/N[a]$ required by Bennett's and Bernstein's inequalities~\citep{AudibertEtAlUCBV}.
This correction term appears to have a significant impact on moderate time horizons: for a sub-optimal arm $a$, the number of draws $N[a]$ does not grow faster than the $\log(t)$ exploration function, and $\log(t)/N[a]$ does not vanish.

\begin{figure}[hbt]
  \centering
  \includegraphics[width=0.96\textwidth]{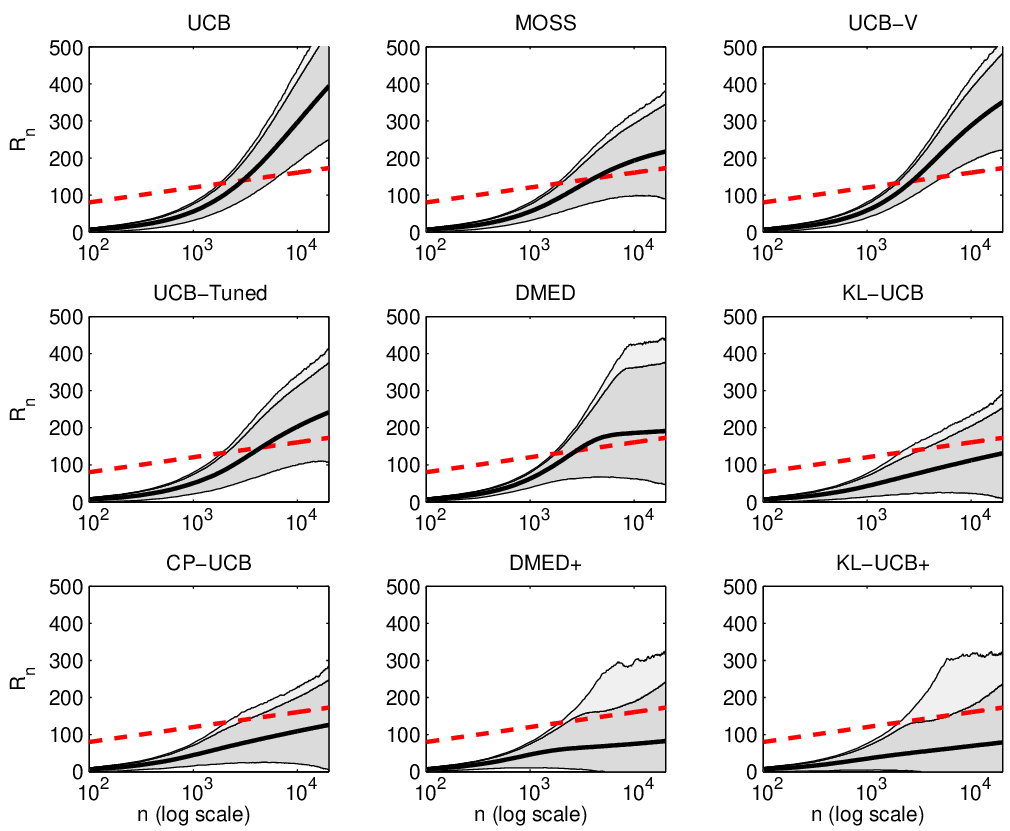}
  \caption{Regret of the various algorithms as a function of time (on a log scale) in the ten arm scenario. On each graph, the red dashed line shows the lower bound, the solid bold curve corresponds to the mean regret while the dark and light shaded regions show respectively the central 99\% region and the upper 0.05\% quantile, respectively.}
  \label{fig:10}
\end{figure}

\subsection{Scenario 2: low rewards}
In Figure~\ref{fig:10} we consider a significantly more difficult scenario,
again with Bernoulli rewards, inspired by a situation (frequent in applications like marketing or Internet advertising) where the mean
reward of each arm is very low. In this scenario, there are ten arms: the
optimal arm has expected reward 0.1, and the nine suboptimal arms consist of
three different groups of three (stochastically) identical arms each with
expected rewards 0.05, 0.02 and 0.01, respectively. We again used $N=50,000$
simulations to obtain the regret plots of Figure~\ref{fig:10}. These plots show, for
each algorithm, the average cumulated regret together with quantiles of the
cumulated regret distribution as a function of time (again on a logarithmic
scale).

In this scenario, the difference is more pronounced between UCB and KL-UCB. The performance gain of
UCB-Tuned is also much less significant. 
KL-UCB and DMED reach a performance
that is on par with the lower bound of \cite{BurnetasKatehakis97bandits}
in~(\ref{eq:binfRegret}), although the performance of KL-UCB is here again
significantly better. Using KL-UCB+ and DMED+ results in significant mean
improvements, although there are hints that those algorithms might indeed be too
risky with occasional very large deviations from the mean regret curve.

The final algorithm included in this roundup, called CP-UCB, is in some
sense a further adaptation of KL-UCB to the specific case of Bernoulli rewards.
For $n\in\N$ and $p\in[0,1]$, denote by $P_{n,p}$
the binomial distribution with parameters $n$ and $p$.  For a random variable
$X$ with distribution $P_{n,p}$, the \emph{Clopper-Pearson}~(see \cite{ClopperPearson34}) or ``exact''
upper-confidence bound of risk $\alpha\in]0,1[$ for $p$ is
\[u^{CP}(X, n, \alpha) = \max\left\{q \in [0,1] : P_{n,q}([0, X])\geq
  \alpha\right\} \;.\] It is easily verified that $P_{n,p}\left(\mu\leq
  u^{CP}(X)\right) \geq 1-\alpha$, and that $u^{CP}(X)$ is the smallest
quantity satisfying this property: $u^{CP}(X)\leq \tilde{u}(X)$ for any other
upper-confidence bound $\tilde{u}(X)$ of risk at most $\alpha$.

The Clopper-Pearson Upper-Confidence Bound algorithm (CP-UCB) differs from KL-UCB only in the way the upper-confidence bound on the performance of each arm is computed, replacing line 6 of Algorithm~\ref{alg:klUCB} by
 \[a \gets \argmax_{1\leq a\leq K}
    \;u^{CP}\left(S[a], N[a], \frac{1}{t\log(t)^c}\right)\;.\]
 As the Clopper-Wilson confidence intervals are always sharper than the Kullback-Leibler
intervals, one can very easily adapt the proof of Section~\ref{sec:proof} to
show that the regret bounds proved for the KL-UCB algorithm also hold for
CP-UCB in the case of Bernoulli rewards. However, the improvement over KL-UCB is very limited (often, the two algorithms actually take exactly the same decisions). In terms of results, one can observe on Figure~\ref{fig:10} that CP-UCB only achieves a performance that is marginally better than that of KL-UCB. Besides, there is no guarantee that the CP-UCB algorithm is also efficient on arbitrary bounded distributions.

\begin{figure}[hbt]
  \centering
  \includegraphics[width=0.60\textwidth]{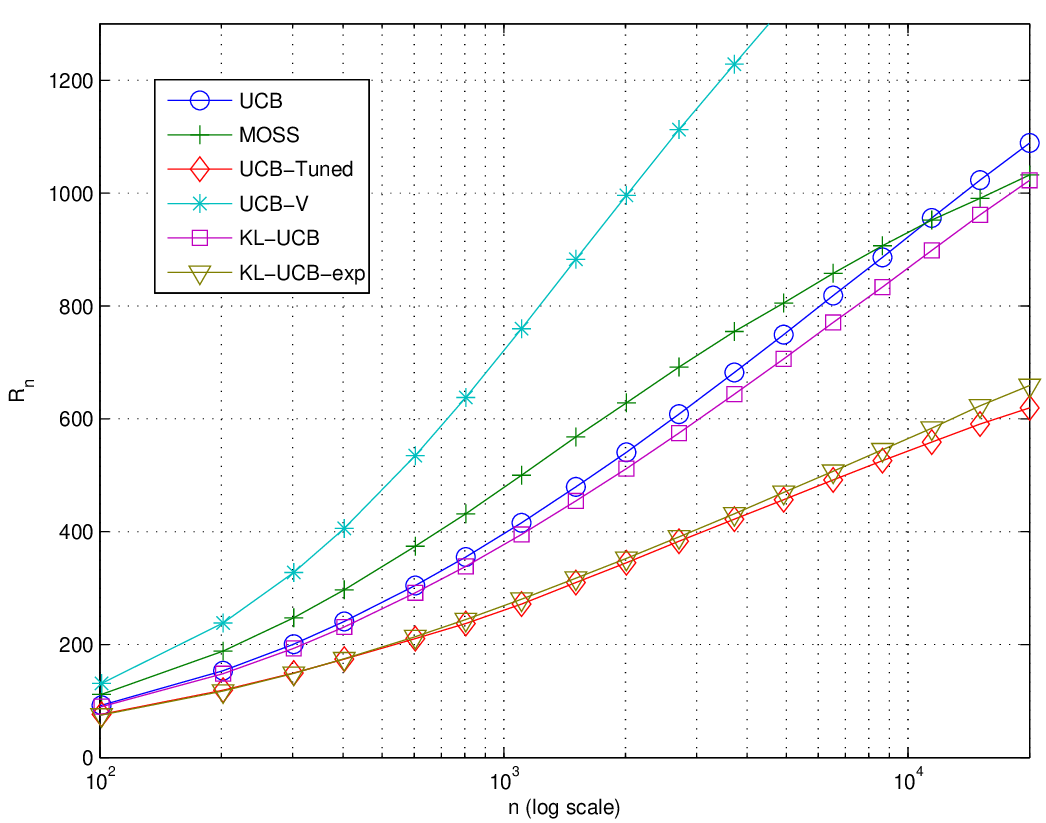}
  \caption{Regret of the various algorithms as a function of time in the bounded exponential scenario.}
  \label{fig:EB5}
\end{figure}

\subsection{Scenario 3: bounded exponential rewards}
In the third example, there are $5$ arms: the rewards are exponential variables, with parameters $1/5$, $1/4$, $1/3$, $1/2$ and $1$ respectively, truncated at $x_{\max}=10$ (thus, they are bounded in $[0,10]$).
The interest of this scenario is twofold: first, it shows the performance of KL-UCB for non-binary, non-discrete, non $[0,1]$-valued rewards. Second, it illustrates that, as explained in Section~\ref{sec:parametricKLUCB}, specific variants of the KL-UCB algorithm can reach an even better performance.

In this scenario, UCB and MOSS, but also KL-UCB are clearly sub-optimal.
UCB-Tuned and UCB-V, by taking into account the variance of the reward distributions (much smaller than the variance of a $\{0,10\}$-valued distribution with the same expectation), were expected to perform significantly better. For the reasons mentionned above this is not the case for UCB-V on a time horizon $n=20,000$. Yet, UCB-Tuned is spectacularly more efficient, and is only caught up by KL-UCB-exp, the variant of KL-UCB designed for exponential rewards. 
Actually, the KL-UCB-exp algorithm ignores the fact that the rewards are truncated, and uses the divergence $d(x,y)=x/y-1-\log(x/y)$ prescribed for genuine exponential distributions. One can easily show that this choice leads to slightly too large upper confidence bounds. Yet, the performance is still excellent, stable, and the algorithm is particularly simple. 

\section{Proof of Theorem~\ref{th:borneNbplayed}}\label{sec:proof}
Consider a positive integer $n$, a small $\epsilon>0$, an optimal arm $a^*$ and
a sub-optimal arm $a$ such that $\mu_a < \mu_{a^*}$.  Without loss of
generality, we will assume that $a^* = 1$. For any arm $b$, the past average performance of arm $b$ is denoted
by $\hat{\mu}_b(t) = S_b(t) / N_b(t)$; by convenience, for every positive
integer $s$ we will also denote $\hat{\mu}_{b,s} = \left( X_{b,1} + \dots +
  X_{b,s} \right)/s$, so that $\hat{\mu}_t(b) = \hat{\mu}_{b,N_b(t)}$. KL-UCB
relies on the following upper-confidence bound for $\mu_b$:
\[u_b(t) = \max \left\{q>\hat{\mu}_b(t) : N_b(t)\,d\left( \hat{\mu}_b(t), q
  \right) \leq \log(t) + 3\log(\log(t))\right\}\;.\]
For $x,y\in[0,1]$, define $d^+(x,y) = d(x,y)\1_{x<y}$.
The expectation of $N_n(a)$ is upper-bounded by using the following
decomposition:
\begin{align*}
  \E[N_n(a)] & = \E\left[\sum_{t=1}^n \1\{A_t = a\}\right]\leq \E\left[\sum_{t=1}^n \1\{\mu_1>u_1(t)\}\right] + \E\left[\sum_{t=1}^n \1\{A_t = a, \mu_1\leq u_1(t)\}\right] \\
  & \leq \sum_{t=1}^n \P\left(\mu_1>u_1(t)\right) + \E\left[ \sum_{s = 1}^n
    \1\{ sd^+(\hat{\mu}_{a,s}, \mu_1) < \log(n) + 3\log(\log(n))\}\right]\;,
\end{align*}
where the last inequality is a consequence of Lemma~\ref{lem:majN}.  The first
summand is upper-bounded as
follows: 
by Theorem~\ref{th:borneDev} (proved in the Appendix),
\begin{align*}
  P(\mu_1>u_1(t))& \leq e\left\lceil \log(t)\left( \log(t) + 3\log(\log(t))
    \right)\right\rceil \exp(-\log(t) - 3\log(\log(t))) \\&= \;
  \frac{e\left\lceil \log(t)^2 + 3\log(t)\log(\log(t))\right\rceil}{t
    \log(t)^3}\;.
\end{align*}
Hence,
\[ \sum_{t=1}^n P(\mu_1>u_1(t)) \leq \sum_{t=1}^n \frac{e\left\lceil \log(t)^2
    + 3\log(t)\log(\log(t)) \right\rceil}{t\log(t)^3} \leq C'_1 \log(\log(n))
\]
for some positive constant $C'_1$ ($C'_1\leq 7$ is sufficient).
For the second summand, let 
 \[K_n = \left\lfloor\frac{1+\epsilon}{d^+(\mu_a, \mu_1)}\Big(  \log(n)+3\log(\log(n)) \Big)\right\rfloor\;.\]
 Then:
\begin{align*}
  \sum_{s = 1}^n \P\big( sd^+(\hat{\mu}_{a,s}, \mu_1) < \log(n) &+ 3\log(\log(n))
  \big)\\
  &\leq K_n + \sum_{s=K_n+1}^\infty \P\Big(  sd^+(\hat{\mu}_{a,s}, \mu_1) < \log(n) + 3\log(\log(n)) \Big) \\
  &\leq K_n + \sum_{s=K_n+1}^\infty \P\Big(  K_nd^+(\hat{\mu}_{a,s}, \mu_1) <\log(n) + 3\log(\log(n)) \Big) \\
  &= K_n + \sum_{s=K_n+1}^\infty \P\left( d^+(\hat{\mu}_{a,s}, \mu_1) < \frac{d(\mu_a, \mu_1)}{1+\epsilon} \right) \\
  &\leq \frac{1+\epsilon}{d^+(\mu_a, \mu_1)}\Big( \log(n)+3\log(\log(n)) \Big)
  + \frac{C_2(\epsilon)}{n^{\beta(\epsilon)}}
\end{align*}
according to Lemma~\ref{lem:majSumProba}.
This will conclude the proof, provided that we prove the following two lemmas.

\begin{lemma}\label{lem:majN}
  \[ \sum_{t=1}^n \1\{A_t = a, \mu_1\leq u_1(t)\} \leq \sum_{s = 1}^n \1\{
  sd^+(\hat{\mu}_{a,s}, \mu_1) < \log(n) + 3\log(\log(n))\}\;.\]
\end{lemma}
\begin{proof}
  Observe that $A_t = a$ and $\mu_1\leq u_1(t)$ together imply that $u_a(t)
  \geq u_1(t) \geq \mu_1$, and hence that 
  \[d^+\left(\hat{\mu}_a(t), \mu_1\right) \leq d(\hat{\mu}_a(t), u_a(t)) =
  \frac{\log(t) + 3\log(\log(t))}{N_a(t)}\;.\] 
Thus,
  \begin{align*}
    \sum_{t=1}^n &\1\{A_t = a, \mu_1\leq u_1(t)\}
     \leq  \sum_{t=1}^n \1\{A_t = a, N_a(t)d^+(\hat{\mu}_a(t), \mu_1) \leq \log(t) + 3\log(\log(t))\}  \\
    & =  \sum_{t=1}^n \sum_{s=1}^t \1\{N_t(a)=s, A_t = a, sd^+(\hat{\mu}_{a,s}, \mu_1) \leq \log(t) + 3\log(\log(t))\} \\
    & \leq  \sum_{t=1}^n \sum_{s=1}^t \1\{N_t(a)=s, A_t = a\}\1\{sd^+(\hat{\mu}_{a,s}, \mu_1) \leq \log(n) + 3\log(\log(n))\} \\
    & = \sum_{s=1}^n \1\{sd^+(\hat{\mu}_{a,s}, \mu_1) \leq \log(n) + 3\log(\log(n))\}  \sum_{t=s}^n \1\{N_t(a)=s, A_t = a\}\\
    & = \sum_{s=1}^n \1\{sd^+(\hat{\mu}_{a,s}, \mu_1) \leq \log(n) +
    3\log(\log(n))\} \;,
  \end{align*}
  as, for every $s\in\{1,\dots,n\}$, $\sum_{t=s}^n \1\{N_t(a)=s, A_t = a\}\leq
  1$.
\end{proof}

\begin{lemma}\label{lem:majSumProba}
  For each $\epsilon>0$, there exist $C_2(\epsilon)>0$ and $\beta(\epsilon)>0$
  such that
  \[\sum_{s=K_n+1}^\infty \P\left(d^+(\hat{\mu}_{a,s}, \mu_1) < \frac{d(\mu_a,
      \mu_1) }{1+\epsilon} \right) \leq
  \frac{C_2(\epsilon)}{n^{\beta(\epsilon)}}\;.\]
\end{lemma}
\begin{proof}
  If $d^+(\hat{\mu}_{a,s}, \mu_1) < d(\mu_a, \mu_1)/(1+\epsilon)\;$, then
  $\hat{\mu}_{a,s}> r(\epsilon)$, where $r(\epsilon)\in ]\mu_a, \mu_1[$ is such
  that $d(r(\epsilon), \mu_1) = d(\mu_a, \mu_1)/(1+\epsilon)$.
  Hence, \begin{multline*}\P\left(d^+(\hat{\mu}_{a,s}, \mu_1)< \frac{d(\mu_a,
        \mu_1) }{1+\epsilon}\right) \leq \P\left( d(\hat{\mu}_{a,s}, \mu_a) >
      d(r(\epsilon), \mu_a),\;\hat{\mu}_{a,s}>\mu_a \right) \\\leq
    \P(\hat{\mu}_{a,s}> r(\epsilon)) \leq \exp(-s d(r(\epsilon), \mu_a))
    \;,\end{multline*} and
  \[\sum_{s=K_n+1}^\infty \P\left(d^+(\hat{\mu}_{a,s}, \mu_1) < \frac{d(\mu_a,
      \mu_1)}{1+\epsilon} \right) \leq \frac{\exp(-d(r(\epsilon),\mu_a)
    K_n)}{1-\exp(-d(r(\epsilon),\mu_a))} \leq
  \frac{C_2(\epsilon)}{n^{\beta(\epsilon)}}\;,\] with
  $C_2(\epsilon)=(1-\exp(-d(r(\epsilon),\mu_a)))^{-1}$ and $\beta(\epsilon) =
  (1+\epsilon)d(r(\epsilon), \mu_a)/d(\mu_a, \mu_1)$.  Easy computations show
  that $r(\epsilon) = \mu_a+O(\epsilon)$, so that $C_2(\epsilon) =
  O(\epsilon^{-2} )$ and $\beta(\epsilon) = O(\epsilon^{2})$.
\end{proof}

\section{Conclusion}\label{sec:conclusion}
The self-normalized deviation bound of Theorems~\ref{th:borneDev} and~\ref{th:borneDevGen}, together with the new analysis presented in Section~\ref{sec:proof}, allowed us to design and analyze improved UCB algorithms.
In this approach, only an upper-bound of the deviations (more precisely, of the exponential moments) of the rewards is required, which makes it possible to obtain versatile policies satisfying interesting regret bounds for large classes of reward distributions.
The resulting index policies are simple, fast, and very efficient in practice, even for small time horizons. 
\bibliography{klUCB}

\appendix
\section{Kullback-Leibler deviations for bounded variables with a random number of summands}\label{sec:appendix}
We start with a simple lemma justifying the focus on binary rewards.
\begin{lemma}\label{lem:bounded2bernoulli}
  Let $X$ be a random variable taking value in $[0,1]$, and let $\mu =
  \E[X]$. Then, for all $\lambda \in \R$, \[E\left[ \exp(\lambda X) \right]
  \leq 1-\mu +\mu\exp(\lambda) \;,\]

\end{lemma}
\begin{proof}
  The function $f:[0,1]\to\R$ defined by
  $f(x) = \exp(\lambda x)-x\left(\exp(\lambda)-1\right) - 1$ is convex and such that $f(0)=f(1)=0$, hence $f(x)\leq 0$ for all
  $x\in[0,1]$.  Consequently,
  \[\E\left[ \exp(\lambda X) \right] \leq \E\left[
    X\left(\exp(\lambda)-1\right) +1 \right] = \mu(\exp(\lambda)-1) +1\;.\]
\end{proof}

\begin{theorem}\label{th:borneDev}
  Let $(X_t)_t\geq 1$ be a sequence of independent random variables bounded in
  $[0,1]$ defined on a probability space $(\Omega, \F, \P)$ with common
  expectation $\mu=\E[X_t]$.  Let $\F_t$ be an increasing sequence of
  $\sigma$-fields of $\F$ such that for each $t$, $\sigma(X_1\ldots,
  X_t)\subset \F_t$ and for $s>t$, $X_s$ is independent from $\F_t$.  Consider
  a previsible sequence $(\epsilon_t)_{t\geq 1}$ of Bernoulli variables (for
  all $t > 0$, $\epsilon_t$ is $\F_{t-1}$-measurable). Let $\delta>0$ and for
  every $t\in\{1,\dots,n\}$ let
  \begin{align*}
    &S(t) =\sum_{s=1}^t \epsilon_s X_s \; , \qquad N(t) = \sum_{s=1}^t
    \epsilon_s\;, \qquad
    \hat{\mu}(t) = \frac{S(t)}{N(t)}\;, \\
    &u(n) = \max \left\{q>\hat{\mu}_n : N(n)d\left( \hat{\mu}(n), q \right)
      \leq \delta\right\}\;.
  \end{align*}
  Then
  \[\P\left( u(n)< \mu\right) \leq e\left\lceil \delta\log(n) \right\rceil
  \exp(-\delta)\;.\]
\end{theorem}
\begin{proof}
  For every $\lambda\in\R$, let
  $\phi_\mu(\lambda) = \log \E\left[ \exp\left( \lambda X_1 \right)
  \right]$. By Lemma~\ref{lem:bounded2bernoulli}, it holds that $\phi_\mu(\lambda)\leq \log \left( 1-\mu+\mu \exp\left( \lambda \right) \right)$.  Let
  $W^{\lambda}_0=1$ and for $t\geq 1$, $$W^{\lambda}_t = \exp(\lambda S(t) -
  N(t) \phi_\mu(\lambda)).$$ $\left(W^{\lambda}_t\right)_{t\geq 0}$ is a
  super-martingale relative to $\left( \F_t \right)_{t\geq 0}$. In fact,
  \begin{multline*}
    \E\big[ \exp\left( \lambda\left\{ S(t+1-S(t) \right\} \right) | \F_t\big] =\E\big[ \exp\left( \lambda \epsilon_{t+1} X_{t+1}\right) | \F_t\big]
     =  \exp\big(  \epsilon_{t+1} \log \E\left[ \exp\left( \lambda X_1 \right) \right] \big)\\
     \leq \exp\big(   \epsilon_{t+1} \phi_\mu\left( \lambda \right) \big)
    = \exp \big( \left\{ N(t+1)-N(t) \right\} \phi_\mu\left( \lambda \right)
    \big)
  \end{multline*}
  which can be rewritten as
  \[
  \E\big[\exp\left( \lambda S(t+1)-N(t+1)\phi_\mu\left( \lambda \right)
    \right) | \F_{t}\big]  \leq \exp\left( \lambda S(t)-N(t)\phi_\mu\left(
      \lambda \right) \right).
  \]

  To proceed, we make use of the so-called 'peeling trick' (see for instance
  \cite{massart2007}): we divide the interval $\{1,\dots, n\}$ of possible
  values for $N(n)$ into "slices" $\{t_{k-1}+1,\dots,t_k\}$ of geometrically
  increasing size, and treat the slices independently.  We may assume that
  $\delta>1$, since otherwise the bound is trivial. Take\footnote{if
    $\delta\leq 1$, it is easy to check that the bound holds whatsoever.} $\eta
  = 1/(\delta-1)$, let $t_0=0$ and for $k\in\N^{*}$, let $t_k = \left\lfloor
    (1+\eta)^k \right\rfloor$.  Let $D$ be the first integer such that $t_D\geq
  n$, that is $D=\left\lceil\frac{\log n}{\log 1+\eta}\right\rceil$.  Let $A_k
  = \left\{t_{k-1} < N(n)\leq t_k\right\}\cap \left\{ u(n)<\mu\right\}$. We
  have:
  \begin{equation}
    \label{eq:peeling:basic}
    \P\left( u(n)<\mu\right)  \leq \P\left( \bigcup_{k=1}^D A_k\right) \leq \sum_{k=1}^D \P\left(A_k\right).
  \end{equation}
  Observe that $u(n)<\mu $ if and only if $\hat{\mu}(n)<\mu$ and
  $N(n)d(\hat{\mu}(n), \mu)>\delta$.  Let $s$ be the smallest integer such that
  $\delta/(s+1)\leq d(0;\mu)$; if $N(n)\leq s$, then $N(n)d(\hat{\mu}(n),
  \mu)\leq s d(\hat{\mu}(n), \mu)\leq sd(0,\mu)< \delta$ and $\P( u(t)< \mu
  )=0$.  Thus, $\P\left(A_k\right)=0$ for all $k$ such that $t_k\leq s$.

  For $k$ such that $t_k>s$, let $\tilde{t}_{k-1} = \max\{t_{k-1}, s\}$.
  Let $x \in ]0, \mu[$ be such that $d(x;\mu)=\delta/N(n)$ and let
  $\lambda(x)=\log(x\left( 1-\mu \right))-\log (\mu\left( 1-x \right))<0$, so
  that $d(x; \mu) =  \lambda(x)x- \left( 1-\mu+\mu \exp\left( \lambda(x) \right) \right).$ Consider $z$ such that $z <
  \mu$ and $d(z, \mu)=\delta/(1+\eta)^{k}$. Observe that:
  \begin{itemize}
  \item if $N(n)>\tilde{t}_{k-1}$, then $$d(z; \mu) =
    \frac{\delta}{(1+\eta)^{k}} \geq \frac{\delta}{(1+\eta)N(n)}\;;$$
  \item if $N(n)\leq t_k$, then as
 $$d(\hat{\mu}(n); \mu)>\frac{\delta}{N(n)}>\frac{\delta}{(1+\eta)^k} = d(z; \mu),$$
it holds that :
  \[\hat{\mu}(n) < \mu \hbox{ and } d(\hat{\mu}(n); \mu)>\frac{\delta}{N(n)}\implies \hat{\mu}(n)\leq z.\]
\end{itemize}
Hence, on the event $\left\{\tilde{t}_{k-1}<N(n)\leq t_k\right\} \cap
\left\{\hat{\mu}(n) < \mu\right\} \cap \left\{d(\hat{\mu}(n);
  \mu)>\frac{\delta}{N(n)}\right\}$ it holds that
$$\lambda(z) \hat{\mu}(n) - \phi_\mu(\lambda(z)) \geq 
\lambda(z) z - \phi_\mu(\lambda(z)) = d(z; \mu)
\geq\frac{\delta}{(1+\eta)N(n)}.$$ Putting everything together,
\begin{align*}
  \left\{ \tilde{t}_{k-1}<N(n)\leq t_k \right\}\cap&\left\{ \hat{\mu}(n) < \mu  \right\} \cap\left\{d(\hat{\mu}(n); \mu) \geq \frac{\delta}{N(n)} \right\}\\
 & \subset \left\{ \lambda(z) \hat{\mu}(n) - \phi_\mu(\lambda(z)) \geq  \frac{\delta}{N(n)\left( 1+\eta \right)}  \right\}\\
   &\subset \left\{  \lambda(z)S_n - N(n) \phi_\mu(\lambda(z)) \geq \frac{\delta}{1+\eta} \right\} \\
   &\subset \left\{ W^{\lambda(z)}_n > \exp
    \left(\frac{\delta}{1+\eta}\right)\right\}.
\end{align*}
As $\left(W^{\lambda}_t\right)_{t\geq0}$ is a supermartingale,
$\E\left[W^{\lambda(z)}_n\right] \leq \E\left[W^{\lambda(z)}_0\right] =1$, and
the Markov inequality yields:
\begin{multline*}
  \P\Big(\left\{\tilde{t}_{k-1}<N(n)\leq t_k\right\}\cap \left\{
      \hat{\mu}(n)\geq \mu \right\} \cap \left\{N(n)d(\hat{\mu}(n), \mu)\geq
      \delta \right\} \Big) \\\leq
  \P\left(W^{\lambda(z)}_n > \exp\left(\frac{\delta}{1+\eta}\right) \right) \label{eq:umarkov} 
  \leq \exp\left(-\frac{\delta}{1+\eta}\right).\nonumber
\end{multline*}
Finally, by Equation \eqref{eq:peeling:basic},
$$ \P\left( \bigcup_{k=1}^D \left\{\tilde{t}_{k-1}<N(n)\leq t_k\right\} \cap \left\{ u(n)<\mu \right\}\right)  \leq D \exp\left(-\frac{\delta}{1+\eta}\right).$$
But as $\eta = 1/(\delta-1)$, $D=\left\lceil\frac{\log n}{\log \left(
      1+1/(\delta-1) \right)}\right\rceil$ and as $\log(1+1/(\delta-1))\geq
1/\delta$, we obtain:
 $$\P\left( u(n)<\mu\right) \leq \left\lceil\frac{\log n}{\log \left( 1+\frac{1}{\delta-1}\right)}\right\rceil 
 \exp(-\delta+1)\leq e\left\lceil \delta\log(n) \right\rceil \exp(-\delta).
 $$ 
\end{proof}

Of course, a symmetric bound for the probability of over-estimating $\mu$ can be
derived following the same lines. Together, they show that for all $\delta>0$:
\[\P\big( N(n)d(\hat{\mu}(n), \mu)>\delta \big) \leq 2e\left\lceil
  \delta\log(n) \right\rceil \exp(-\delta)\;.\]
Finally, we state a more general deviation bound for arbitrary reward distributions with finite exponential moments. The proof (very similar to that of Theorem~\ref{th:borneDev}) is omitted.
\begin{theorem}\label{th:borneDevGen}
  Let $(X_t)_t\geq 1$ be a sequence of i.i.d. random variables defined on a
  probability space $(\Omega, \F, \P)$ with common expectation $\mu$.
  Assume that the cumulant-generating function
$$\phi(\lambda)=\log \E\left[\exp(\lambda X_1)\right]$$
is defined and finite on some open subset $]\lambda_1, \lambda_2[$ of $\R$
containing $0$.   Define $d: \R\times\R \to
\R\cup\{+\infty\}$ as follows:
for all $x\in\R,$
\begin{equation*}
 d(x, \mu) = \sup_{\lambda \in ]\lambda_1, \lambda_2[} \left\{ \lambda x - \phi(\lambda) \right\} \;.
\end{equation*}
Let $\F_t$ be an increasing sequence of $\sigma$-fields of $\F$ such that for
each $t$, $\sigma(X_1\ldots, X_t)\subset \F_t$ and for $s>t$, $X_s$ is
independent from $\F_t$.  Consider a previsible sequence $(\epsilon_t)_{t\geq
  1}$ of Bernoulli variables (for all $t > 0$, $\epsilon_t$ is
$\F_{t-1}$-measurable). Let $\delta>0$ and for every $t\in\{1,\dots,n\}$ let
\begin{align*}
  &S(t) =\sum_{s=1}^t \epsilon_s X_s \; , \qquad N(t) = \sum_{s=1}^t
  \epsilon_s\;, \qquad
  \hat{\mu}(t) = \frac{S(t)}{N(t)}\;, \\
  &u(n) = \max \left\{q>\hat{\mu}_n : N(n)d\left( \hat{\mu}(n), q \right)
    \leq \delta\right\}\;.
\end{align*}
Then
\[\P\left( u(n)< \mu\right) \leq e\left\lceil \delta\log(n)
\right\rceil\exp(-\delta)\;.\]
\end{theorem}

\end{document}